\documentclass[11pt,amstex]{article}
\usepackage{svg}
\usepackage{amscd}
\usepackage{amsmath}
\usepackage{amssymb}
\usepackage{graphicx}

\newtheorem {theo} {\bf Theorem} [section]
\newtheorem {prop} [theo] {\bf Proposition}
\newtheorem {cory} [theo] {\bf Corollary}
\newtheorem {lem} [theo] {\bf Lemma}
\newtheorem {defn} [theo] {\bf Definition}

\newtheorem {rem} [theo] {\bf Remark}
\newcommand{\QED}{\hfill \CaixaPreta \vspace{6mm}}
\def\CaixaPreta{\vrule Depth0pt height6pt width6pt}

\newcommand{\qed}{\nopagebreak\hfill{\vrule width6pt height6pt depth0pt}}
\newcommand{\be}{\begin{eqnarray}}
\newcommand{\ee}{\end{eqnarray}}
\newcommand{\benn}{\begin{eqnarray*}}
\newcommand{\eenn}{\end{eqnarray*}}
\newcommand{\bse}{\begin{equation}}
\newcommand{\ese}{\end{equation}}
\newcommand{\bsenn}{\begin{displaymath}}
\newcommand{\esenn}{\end{displaymath}}
\newcommand{\logand}{\;\;{\rm and }\;\;}

\newcommand{\logif}{\;\;{\rm if }\;\;}

\newcommand{\N}{\mathbb{N}}

\newcommand{\R}{\mathbb{R}}

\newcommand{\Z}{\mathbb{Z}}
\usepackage{authblk}
\begin{document}

\title{Birkhoff Measures, Birkhoff Sums, and Discrepancies}
%\author{J. J. P. Veerman, D. Ralston, F. M. Tangerman, H. Wu}
\author[1]{D. Ralston}
\author[2]{F.M. Tangerman}
\author[3]{J.J.P. Veerman}
\author[4]{H. Wu}
\affil[1]{ralstond@oldwestbury.edu, Math, SUNY Old Westbury, Old Westbury, NY, USA.}
\affil[2]{fmtangerman@gmail.com, FM Tangerman LLC, Hartford CT, USA.}
\affil[3]{veerman@pdx.edu, Math, Portland State University, Portland, OR, USA.}
\affil[4]{haowu.nankai@gmail.com, Inst. of Analysis and Number Theory,
Graz UT, AUSTRIA.}

\maketitle

 %\normalsize        %\mysetfontsize5

\noindent
\section*{Abstract}
\begin{small} We study the distribution of a sequence of points in the circle generated
by rotations by a fixed  irrational number $\rho$ with initial condition $x_0$, that is:
$\{x_0+i\rho\}_{i=1}^n$. The \emph{discrepancy} as defined
by Pisot and Van Der Corput \cite{VdCP}, quantifies how evenly distributed such a sequence is.

Consider the ergodic or Birkhoff sum of mean zero $S(\rho,n,x):=\sum_{i=1}^{n} (\{x+i\rho\}-1/2)$, where
$\{\cdot\}$ denotes the fractional part. This is a piecewise-linear map in the variable $x$ with $n$
branches, each with slope $n$. For fixed $n$ and $\rho$, let $\nu(\rho,n,z)$ be the number of pre-images
of $S(\rho,n,x)=z$ divided by $n$. Then $\nu(\rho,n,z)$ is a probability density. We call the associated
measures Birkhoff measures.

We prove that the length of the support of the Birkhoff measure $\nu(\rho,n,z)dz$ can be expressed in
terms of the discrepancy. We also show that if $n$ is a continued fraction denominator of $\rho$,
then the graph of $\nu(\rho,n,z)$ is an approximate isosceles trapezoid.
We also give new proofs of two classical results, one by Ramshaw
\cite{Ramshaw} and one by Kuipers-Niederreiter \cite{KN}. These results
allow efficient computation of both Birkhoff sums and discrepancies.
\end{small}

%\vskip 0.4in\noindent
\begin{centering}
\section{Introduction}
\label{chap:intro}
\end{centering}
\setcounter{figure}{0} \setcounter{equation}{0}

This paper considers
partial sums of the function $\{x\}-1/2$ composed of a rotation by angle $\rho$, where $\{\cdot\}$ denotes the fractional part. These are also referred to as \emph{ergodic sums} or \emph{Birkhoff sums}
\bse
S(\rho,n,x) := \sum_{i=1}^{n} (\{x+i\rho\}-1/2)=nx+\dfrac{n(n+1)}{2}\,\rho - \dfrac n2 - \sum_{i=1}^{n}\,\lfloor x+i\rho \rfloor \,.
\label{eq:sumrotations}
\ese
We note that in the definition of $S$ the summation starts at 1 (and not at 0) and ends at $n$.
When $\rho$ and $n$ are fixed we may abbreviate $S(\rho,n,x)$ by $S(x)$. For much of this paper
the standing assumption is that $\rho$ is irrational. We will always assume that $\rho\in[0,1)$.

Discrepancy was originally introduced by Pisot and Van Der Corput \cite{VdCP} and later in
\cite{KN}, as a way to measure the extent to which a finite set of points in a given interval is equally
distributed. Let $\overline{x}:=\{x_i\}_{i=1}^\infty$ be a sequence of points in the unit interval
$[0,1)$. The subintervals  $I$ of $[0,1)$ we consider are half open and half closed, of the form
$I=[a,b)$ or its complement (recall that $[0,1)$ is our model for the unit circle). Then denote
\bsenn
A(I,n):= \#\{ 1\le i\le n \ | \ x_i \in I\} \,.
\esenn

\vskip -.0in\noindent
\begin{defn} The discrepancy of the first $n$ points of $\overline{x}$ is then
\bsenn
D_n(\overline{x}):=\sup_{I=[a,b)\subseteq [0,1)} \left(\frac{A(I,n)}{n}-\ell(I)\right) \,,
\esenn
where $\ell(I)$ is the length of the interval $I$. (In practice the quantity $nD_n(\overline{x})$
is more convenient.)
\label{def:discrepancy}
\end{defn}

\vskip -.0in\noindent
The discrepancy\footnote{We remark that in most texts $D_n$ is defined
as $\sup_{I\subseteq [0,1)} \left|\frac{A(I,n)}{n}-\ell(I)\right|$
(i.e. absolute value is taken). However, if $I$ and $J$ are complements in $[0,1)$, then
\bsenn
A(I,n)-n\ell(I)+A(J,n)-n\ell(J)=0 \quad \Longrightarrow \quad A(I,n)-n\ell(I)= n\ell(J) - A(J,n)\,.
\esenn
and so we can omit the absolute signs without loss of generality.} measures the evenness of the
distribution; it is large when there are underpopulated or overpopulated intervals. For example, if
$(x_1,\cdots,x_n)=(\frac 1n, \frac 2n,\cdots, \frac nn)$, then $D_n(\overline{x})=\frac{1}{n}$. At
the other extreme, if $(x_1,\cdots,x_n)=(x_1,x_1,\cdots, x_1)$, then $D_n(\overline{x})=1$.
It is clear that the most evenly distributed set of $n$ points corresponds to the first example 
(up to reordering) with discrepancy equal to $\frac{1}{n}$.

\begin{figure}[!ht]
\centering
\includegraphics[width=4.5in]{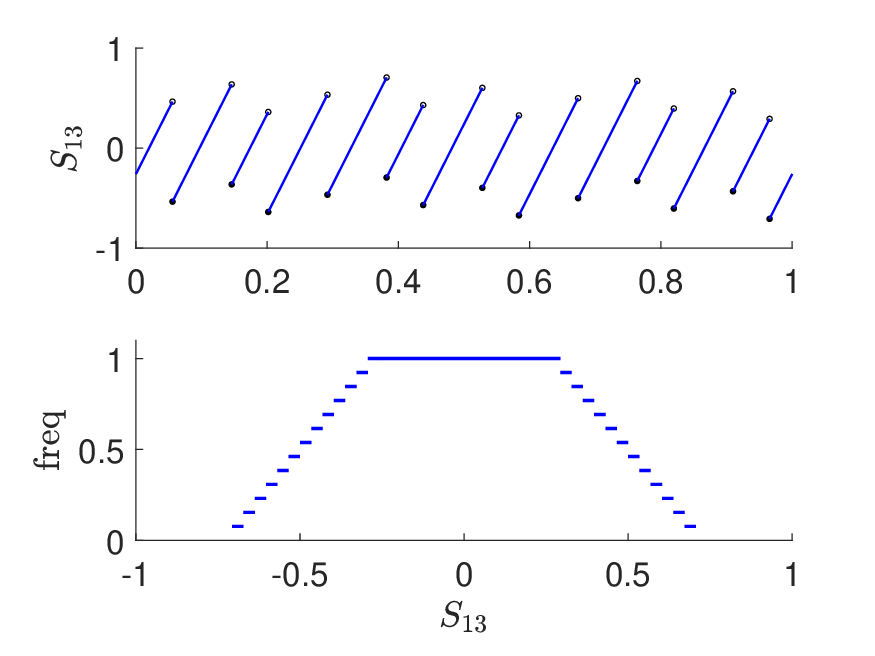}
\caption{\emph{Top figure: $S(\rho,13,x)$ where $\rho$ is the golden mean. Bottom: $\nu(\rho,13,z)$.
Note that $z=0$ intersects all branches, hence $\nu(\rho,13,0)=1$. }}
\label{fig:construction-bhoffmeas}
\end{figure}

One of the aims of this work is to introduce the notion of Birkhoff measures and indicate their
relation with Birkhoff sums and discrepancy. The construction of the associated density is
illustrated in Figure \ref{fig:construction-bhoffmeas}.

\vskip -.0in\noindent
\begin{defn} A Birkhoff density is the pushforward of the Lebesgue measure by $S$:
\bsenn
\nu(\rho,n,z)=\sum_{x:S(\rho,n,x)=z} \frac{1}{|\partial_xS(\rho,n,x)|} \,.
\esenn
Equivalently, it is the number of preimages under $S$ divided by $n$. The associated
measure is the Birkhoff measure.
\label{def:density}
\end{defn}

\vskip -.1in\noindent

\begin{figure}[!ht]
\centering
\includegraphics[width=4.3in]{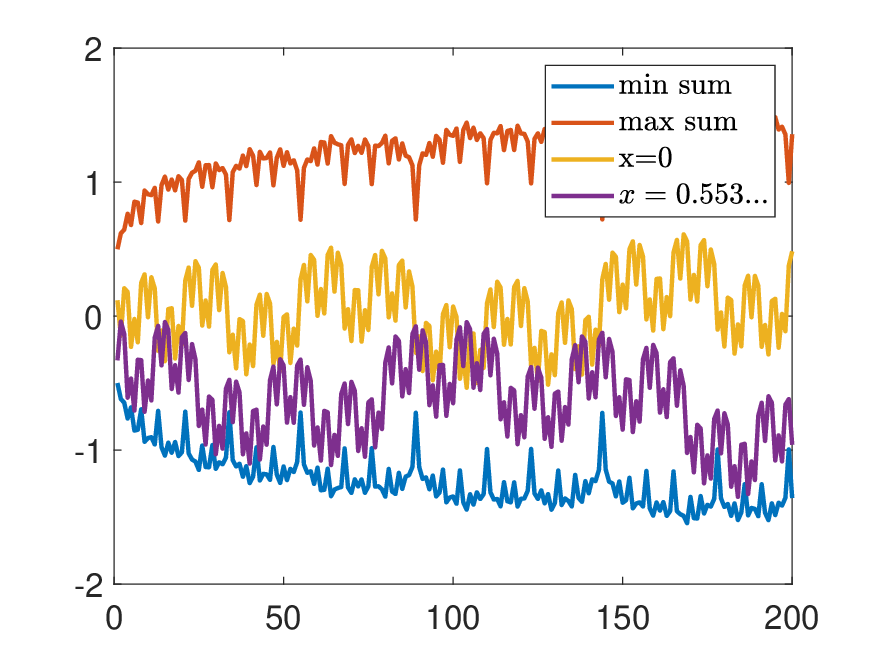}
\caption{\emph{For $\rho$ the golden mean, we show the range of $S(\rho,n,x)$, indicated by the min/max
values for the sum $S$ over $x\in [0,1)$ with increasing value of $n$. The plot also shows
$S(\rho,n,x_0)$ for two initial conditions indicated, $x_0=0$ and $x_0=1-1/\sqrt{5}$.}}
\label{fig:range+orbit}
\end{figure}

\vskip -.1in
We now give a short description of our most important results.
In Figure \ref{fig:range+orbit}, the behavior of $S(\rho, n, x)$ when $\rho$ is the golden mean is
illustrated. The figure shows that the range of $x\rightarrow S(\rho,n,x)$ (which equals the support
of $z\rightarrow \nu(\rho,n,z)$) is symmetric about $0$. In fact,
in Section \ref{chap:measures}, we prove a stronger result, namely that the Birkhoff densities
are symmetric around zero and \emph{tile} the line: they are strictly positive on an interval
symmetric around zero and satisfy $\sum_{m\in\Z}\nu(\rho,n,z+m)=1$ (see Theorem \ref{thm:tile}).
The main result in section \ref{chap:clumpi} is that the length of the support of the Birkhoff
density $\nu(\rho,n,z)$ equals $nD_n(\overline{x})$, where
$\overline{x}=\{i\rho\}_{i=1}^\infty$ (Theorem \ref{thm:clumps=rangeS}).
In section \ref{chap:trapezoid}, we show (Theorem \ref{thm:trapezoid}) that if $p_n/q_n$ is a
continued fraction convergent of $\rho$, then the graph of $\nu(\rho,q_n,z)$ is an approximate
isosceles trapezoid, where the slanted sides consist of a union of small horizontal segments,
as illustrated in the bottom of Figure \ref{fig:construction-bhoffmeas}.

In the appendices \ref{chap:Clump+Bhoff} and \ref{chap:Bhoffrecursion}, we review two important
classical results that provide powerful tools to compute exact values of both the Birkhoff sums and
the discrepancy of sequences generated by rotations. In the first of these, we show how
the length of the support of $S(\rho,n,x)$ can be computed from the sequence $S(\rho,i,0)$ for
$i\in\{1,\cdots,n\}$ (Theorem \ref{thm:subtractingsequences}). Our proof considerably simplifies
the original proof \cite{Ramshaw}. In the latter, our main contribution is the recursive identity in Proposition \ref{prop:recursionS} from which we derive a computational tool (Theorem
\ref{thm:dependence-bk}) that can be used to analyze the growth of $S(\rho,n,0)$ as function of $n$.
This result is due to \cite{KN}; our proof is a more succinct version of the one in \cite{DT}.
In this appendix we also briefly discuss applications of this result. Finally, in appendix
\ref{chap:Cbestiary} we exhibit a sampling of Birkhoff measures with $\rho=e-2$ for various $n$.

\begin{figure}[!ht]
\centering
\includegraphics[height=3.5in,width=7.0in]{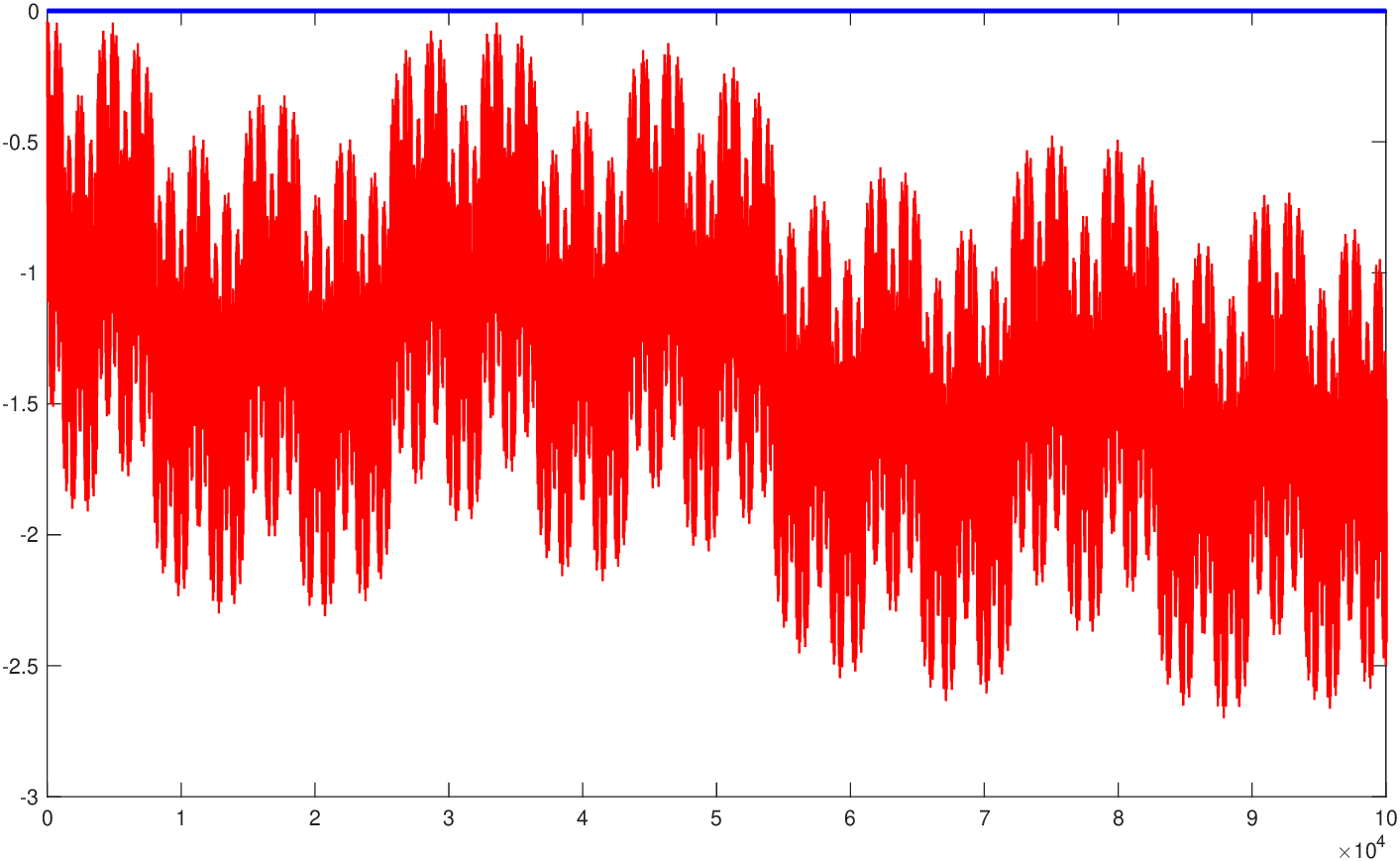}
\caption{\emph{$S(\rho,n,1-1/\sqrt{5})$ where $\rho$ is the golden mean for $n$ up to $10^5$.}}
\label{fig:Snegative}
\end{figure}

To end, we address some open questions. It appears that for certain initial conditions $x$, the
sequence $S(\rho,n,x)$ as function of $n$ is surprisingly \emph{asymmetric}. When $\rho$ is the golden
mean and $x=1-1/\sqrt{5}$, the sequence is \emph{negative} for all $n\leq 10^5$ as we checked in double
precision. We plot the sequence in Figure \ref{fig:Snegative}. Open questions are:
prove that the sequence is negative for all $n$. What are the values $x$ for which
$S$ stays negative?

\vskip -0.2in\noindent

It is known \cite{Aardenne, Aardenne2, Nieder92} that if
$c= \max_{x>3}\frac{x-2}{4(x-1)\ln x}\approx
0.120$, then for any $\overline{x}$ there are infinitely many $n$ such that
$D_n(\overline{x})>c\frac{\ln n}{n}$. Curiously, it
is not clear if there exists a sequence which minimizes that growth rate.

Finally, it would be interesting to understand the shapes of the Birkhoff measures for more values
of $n$ than just the continued fraction approximants. Promising values for $n$ are simple
sums of two continued fraction approximants: $n=2q_k$ or $n=q_k+q_\ell$.

\vskip .0in \noindent
{\bf Acknowledgements:} The authors are grateful to H. Kravitz and H. Moore for developing the
matlab programs used to generate Figures \ref{fig:construction-bhoffmeas}, \ref{fig:range+orbit},
\ref{fig:tiling}, and \ref{fig:reducedresidue}. More details about the algorithm used can be found in \cite{Krav}.

\begin{centering}
\section{Birkhoff Measures and Their Properties}
\label{chap:measures}
\end{centering}
\setcounter{figure}{0} \setcounter{equation}{0}

\begin{prop} The density $\nu(\rho,n,z)$ has the following properties.\\
i) $\nu(\rho,n,z) = \nu(\rho,n,-z) = \nu(1-\rho,n,z)$.\\
ii) The support of $\nu(\rho,n,z)$ is symmetric.\\
iii) $0\leq \nu(\rho,n,z)\leq 1$.
\label{prop:symmetry}
\end{prop}

\vskip -0.0in\noindent
{\bf Proof.} Translation by $-\rho$ starting at $-x$ yields the mirror image of translation by $\rho$
starting at $x$. Thus
\bse
S(1-\rho,n,-x) = S(-\rho,n,1-x) = -S(\rho,n,x) .
\label{eq:symm1}
\ese
This proves that $\nu(1-\rho,n,-z) = \nu(\rho,n,z)$.
It is sufficient to show that $\nu(-\rho,n,z) = \nu(\rho,n,z)$.
From \eqref{eq:sumrotations}, it follows that
\bse
S(1-\rho,n,\{x+{(n+1)}\rho\})=S(-\rho,n,\{x+{(n+1)}\rho\})=S(\rho,n,x) \,,
\label{eq:symm2}
\ese
and so we see that the number of pre-images of $S(\rho,n,\cdot)=z$ equals that of $S(1-\rho,n,\cdot)=z$.
Now the desired symmetry follows from Definition \ref{def:density}. This proves (i).

Item (ii) follows from the first equality of (i). Item (iii) follows from
the fact that $S(\rho,n,x)$ has at most $n$ branches of slope $n$.
\QED

\vskip -0.2in\
The following proposition asserts that $\nu$ is a \emph{tile}: its translates over $\Z$ give the
Lebesgue measure, and $\nu$ itself is (strictly) positive
on an interval and zero elsewhere. These properties are illustrated in
Figure \ref{fig:tiling}.

\begin{figure}[!ht]
\centering
\includegraphics[width=4.5in]{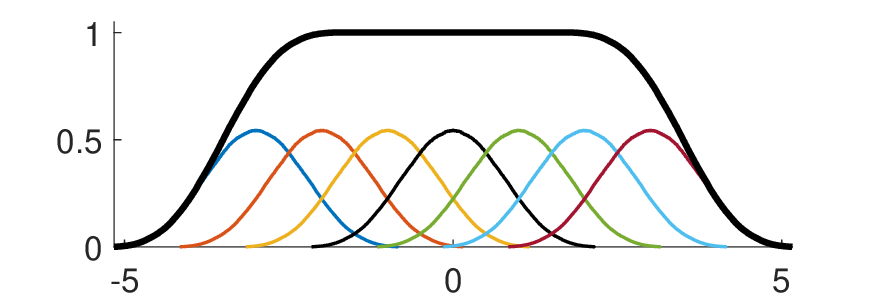}
\caption{\emph{The sum of seven translates of the Birkhoff measure associated with $e-2$ and $n=2024$.}}
\label{fig:tiling}
\end{figure}

\begin{theo}[Tiling Property] Let $\rho$ be irrational and $n$ a positive integer $n$. Then\\
i)  For every $z\in [0,1)$: $\sum_{i\in \Z} \nu(\rho,n,z+i)=1$,\\
ii) The support of the measure associated with the density $\nu(\rho,n,z)$ is a bounded interval.
\label{thm:tile}
\end{theo}

\vskip -0.0in\noindent
{\bf Proof.} Denote the number on pre-images of $z$ under $S(\rho,n,x)$ by $\# S^{-1}(z)$. To prove (i), note that
\bsenn
\begin{aligned}
 & \sum_{i \in \mathbb{Z}} \nu(\rho, n, z+i)=1 \quad \Longleftrightarrow \quad
\sum_{i \in \mathbb{Z}} \# S^{-1}(z+i)=n .
\end{aligned}
\esenn
The right hand side says that we consider pre-images of $\pi\circ S$, where $\pi$ is the canonical
projection of $\R$ to the circle. The branches of $S$ have slope $n$ while at the discontinuities,
$S$ changes by an integer. Therefore, $\pi\circ S= (c+nx) \mod 1$ for some constant $c$. It follows
that there are precisely $n$ inverse images of $\pi \circ S$ for each point on the circle.

To prove (ii), observe that $S$ has the following two properties:\\
(a) its only discontinuities points are located at $\{x_i\}_{i=1}^k$ such
that $\lim_{x\to x_i^-}S(x)>\lim_{x\to x_i^+} S(x)$ (that is: all discontinuities are `downward'
jumps).\\
(b) $S(1^{-}):=\lim_{x\rightarrow 1^{-}}S(x)\geq S(0)$. (Inequality only if $0$ is a discontinuity, in which case $S(1^{-})=S(0)+1$).\\
We prove that the image of $S$ is an interval.

Let $y_+$ ($y_-$) be the supremum (infimum) of the image of $S$, and set $x_\pm$ such that
$\limsup_{x\rightarrow x_+}S(x)=y_+$, and similar for $x_-$. First suppose that 0 is not a discontinuity.
Then 0 is neither a minimum nor a maximum, and so $y_+>S(0)$ and $y_-<S(1^{-})$.
On $[0,x_+]$ define $F(x):=\sup_{0\leq x'\leq x}S(x')$. $F$ is continuous on this interval and so
the image of $F$ equals $[S(0),y_+)$. A fortiori, the image of $S$ contains $[S(0),y_+)$.
By the same reasoning, the image of $S$ must contain all of $(y_-,S(1^{-}))$.
Since $S(1^{-})\geq S(0)$, the image of $S$ must contain all of $(y_-,y_+)$.
The same reasoning work when $S(1^{-})=S(0)+1$. This proves (ii).
\QED

\vskip -0.1in
We know that for any fixed irrational $\rho$, that as $n\rightarrow \infty$, the shape of the graph of
$z\rightarrow \nu(\rho,n,z)$ does
\emph{not} converge. We will see that on the one hand a trapezoid with support in a bounded interval
occurs infinitely often (Theorem \ref{thm:trapezoid}) as shape of $\nu(z)$, while on the
other hand there is also a subsequence $\{n_i\}$ where the corresponding Birkhoff sums $S(\rho, n, x)$,
after normalization by a constant $B_n$, converge to a Gaussian distribution \cite{conze}.

\begin{prop} Let $n$ be fixed. The function from $[0,1)$ to $L^1(\R)$ given by
$\rho\rightarrow \nu(\rho,n,\cdot)$ is continuous:
\bsenn
\lim_{\rho\rightarrow \rho_0}\int_\R|\nu(\rho,n,z) - \nu(\rho_0,n,z)|\,dz = 0 \,.
\esenn
\label{prop:conts-inL1}
\end{prop}

\vskip -0.2in\noindent
{\bf Proof.}  Denote by $D$ the set of points of discontinuity of $S$, consisting of the points
$\{-i\rho\}_{i=1}^n$. For $\epsilon>0$, consider the set $D_{\epsilon}$ of points within $\epsilon/2$ of
$D$. $D_{\epsilon}$ has Lebesgue measure bounded by $n\epsilon$. Since $S$ has derivative $n$ (where
it exists), there is a region $J_{\rho}$ of measure at most $n^2\epsilon$ in $\R$ where
there can be a different number of pre-images under $S(\rho,n,x)$ and $S(\rho_0,n,x)$.
And so for $\rho$ sufficiently close to $\rho_0$:
\benn
\int_{\R}|\nu(\rho,n,z) - \nu(\rho_0,n,z)| dz &= &
 \int_{J_{\rho}}|\nu(\rho,n,z)-\nu(\rho_0,n,z)| dz +
 \int_{\R-J_{\rho}}|\nu(\rho,n,z)-\nu(\rho_0,n,z)|dz\\
&\leq & 2 \int_{J_{\rho}}1 \,dz+\int_{\R-J_{\rho}}0\, dz\leq
2n^2\epsilon\,,
\eenn
where we have used Proposition \ref{prop:symmetry} (iii).
\QED

%\newpage
%\vskip 0.4in\noindent
\begin{centering}
\section{Discrepancy and the Range of Birkhoff Measures }
\label{chap:clumpi}
\end{centering}
\setcounter{figure}{0} \setcounter{equation}{0}

\vskip -0.0in\noindent
\begin{theo} The range of $S(x)$ over $x\in [0,1)$ is an interval of length $nD_n(\overline{x})$, where
$\overline{x}=\{i\rho\}_{i=1}^\infty$.
\label{thm:clumps=rangeS}
\end{theo}

\vskip 0.0in\noindent
{\bf Proof.} Let $-\overline{x}=\{-i\rho\}_{i=1}^\infty$. By symmetry,
$nD_n(\overline{x})=nD_n(-\overline{x})$.
Denote the points $\{-i\rho\}_{i=1}^n$ in $[0,1)$ by $y_i$ in ascending order so that
\bsenn
0\leq y_1\leq y_2 \cdots \leq y_n< 1 \,.
\esenn
Because the definition of $D_n$ involves finitely many points $y_i$, its determination
is then a maximum over finitely many quantities, namely
\bsenn
nD_n(\overline{x})=nD_n(-\overline{x}):=\max_{1< i,j\leq n} \left((j-i+1)-n(y_j-y_i)\right)=
1+\max_{1< i,j\leq n} \left((j-i)-n(y_j-y_i)\right)\,.
\esenn
Definition \ref{def:discrepancy} gives
\bse
nD_n(\overline{x})=1+\max_{0< i,j\leq n} \left(n(y_j-y_i)-(j-i)\right)\,.
\label{eq:clumpiness}
\ese
\noindent
Recall that $S(\rho,n,x)$ has $n$ branches with discontinuities at $\{-i\rho\}_{i=1}^n$. Since each
branch has slope $n$, and is reduced at every discontinuity by one, the difference of the values of
the suprema $S^+(y_i)$ of the successive branches can be computed (for fixed $\rho$ and $n$) as
\bsenn
S^+(y_{i+1})=S^+(y_i)-1+n(y_{i+1}-y_i) =S^+ (y_1)+ n(y_{i+1}-y_1)-i\,.
\esenn
The local minima of $S^-(y_{j+1})$ are exactly one less than the local suprema, so
\bsenn
S^-(y_{j+1})=S^+(y_{j+1})-1 = S^+(y_1)-1+ n(y_{j+1}-y_1)-j\,.
\esenn
The range of $S(x)$ is the interval bounded by the largest of the former and the smallest of the latter.
Its length is the difference of those two, which together with \eqref{eq:clumpiness} yields the theorem.
\QED

\vskip -0.2in\noindent
\begin{cory} The support of $\nu(\rho, n, z)$, as a function of $z$, equals the range of $S(x)$ and
is an interval symmetric around 0 whose length is
equal to $nD_n(\overline{x}_\rho)$, where $\overline{x}_\rho=\{i\rho\}_{i=1}^\infty$.
\label{cor:Cn-is-length-of-support of nu}
\end{cory}

\vskip 0.0in\noindent
{\bf Proof.} The support of $\nu$ equals the range of $S$. The symmetry was established in Proposition
\ref{prop:symmetry} (ii). \QED

\vskip -0.1in\noindent
\begin{cory} For fixed n, the function from $[0,1)$ to $\R$ given by
$\rho\rightarrow nD_n(\overline{x}_\rho)$
where $\overline{x}_\rho=\{i\rho\}_{i=1}^\infty$, is continuous.
\label{cor:Cn-is-conts}
\end{cory}

\vskip 0.0in\noindent
{\bf Proof.} Recall that $\nu(\rho,n,z)\geq 1/n$ on its support. Thus if the support of $\nu(\rho,n,z)$
differs by at least $\epsilon>0$ from that of $\nu(\rho_0,n,z)$ as $\rho\rightarrow \rho_0$, then
the integral in Proposition \ref{prop:conts-inL1} would yield at least $\epsilon/n$.
\QED

%\vskip 0.4in\noindent
\begin{centering}
\section{The Trapezoid Theorem}
\label{chap:trapezoid}
\end{centering}
\setcounter{figure}{0} \setcounter{equation}{0}

\vskip -0.0in\noindent
In this section we show that for continued fraction approximants $q_n$ to $\rho$ that the graph of the
density $\nu(\rho,q_n,.)$ is \textit{approximately trapezoidal}, with the precise definition given
in Definition \ref{def:trapezoid}. First, we state a result and a few conventions that we will need
later. For the conventions we follow \cite{book}, Chapter 6. Recall that we assume $\rho$ to be
irrational.

\vskip -0.0in\noindent
\begin{prop} Let $\gcd(p,q)=1$ and $d=q\rho-p$ and suppose that $|d|<1/(q-1)$. Then
\benn
i)\hskip .3in \sum_{i=1}^{q}\,\left\lfloor i\rho\right\rfloor&=&\frac{(q+1)p-q+1}{2}+\lfloor d \rfloor \\
ii) \hskip .3in \sum_{i=1}^{q}\,\left\{ i\rho\right\}&=&\frac{(q+1)d+q-1}{2}-\lfloor d \rfloor \,.
\eenn
\label{prop:funnysum}
\end{prop}

\vskip -0.3in\noindent
{\bf Proof.} First, set $d=0$ or $\rho=p/q$. Then $\sum_{i=1}^{q}\,\left\{i\rho\right\}$ can be
computed in two ways. The first is $-1+\sum_{i=1}^{q}\, i/q$ (the "-1" term cancels the $q/q$ term
in the sum), the second is $\sum_{i=1}^{q}\,ip/q-\left\lfloor ip/q\right\rfloor$. Equating the two gives
\bsenn
\sum_{i=1}^q\frac iq-1= \frac{q(q+1)}{2}\,\frac pq -
\sum_{i=1}^q\left\lfloor i\frac pq\right\rfloor \,.
\esenn
This yields the expression in item (i) for $d=0$. Item (ii) is then found as the difference
$\sum_{i=1}^q i\rho-\sum_{i=1}^q \lfloor i\rho\rfloor$.

As one varies $\rho$, $\lfloor i\rho\rfloor$ for $i\in \{1,\cdots,q-1\}$ in item (i) can only change if
it passes through an integer value. Since $i\,\frac pq$ has distance at least $1/q$ to any integer for
$1\le i\le q-1$, this means that $\left|i\left(\rho-\frac pq\right)\right|=\left|i\frac dq\right|$
must be at least $1/q$ for the sum to change, which is prevented by the condition on $d$.
The only other possibility is for $\lfloor q\rho\rfloor$ to change, and the change
is equal to $\lfloor q\rho\rfloor-p =\lfloor d\rfloor$. So we must add that term to (i).
\QED

\vskip -0.2in\noindent
\begin{defn} The continued fraction coefficients of $\rho\in(0,1)$
are denoted by $a_n$ so that $\rho=[a_1,a_2,\cdots ]$. The continued fraction convergents to $\rho$ are
denoted by $p_n/q_n$ with $p_0=0$, $p_1=1$, $q_0=1$, and $q_1=a_1$, and for $n\geq 2$,
\bsenn
\begin{array}{ccc}
p_{n+1} &=& a_{n+1}p_{n}+p_{n-1}\\[-0.0cm]
q_{n+1} &=& a_{n+1}q_{n}+q_{n-1} \;.
\end{array}
\esenn
\label{def:contdfr}
\end{defn}

\vskip -0.3in\noindent
\begin{defn} For given $\rho$ and $n\geq 0$, define $d_n$ as
\bsenn
d_n:= q_n\rho-p_n  \,,
\esenn
Note that for $n$ odd, the entry $d_n$ is negative while for $n$ even $d_n$ is positive. Furthermore,
$d_0=\rho$, $d_1=a_1\rho-1$, and for $n\geq 2$
\bsenn
\begin{array}{ccc}
d_{n+1} &=& a_{n+1}d_{n}+d_{n-1}
\end{array} \;.
\esenn
\label{def:dn}
\end{defn}

\vskip -0.3in\noindent
\begin{prop} Let $d=q\rho-p$ and suppose that $|d|<1/(q-1)$, with $p$ and $q$ coprime. We have:
\bsenn
\nu(\rho,q,z) = 1 \quad \textrm{iff}\quad z\in \left[\frac{-(1-(q-1)|d|)}{2},\frac{(1-(q-1)|d|)}{2}\right)\,,
\esenn
and $\nu(\rho,q,z)\in [0,\frac{q-1}{q}]$ otherwise.
\label{prop:nu=1}
\end{prop}

\begin{figure}[!ht]
\centering
\includegraphics[width=4.0in]{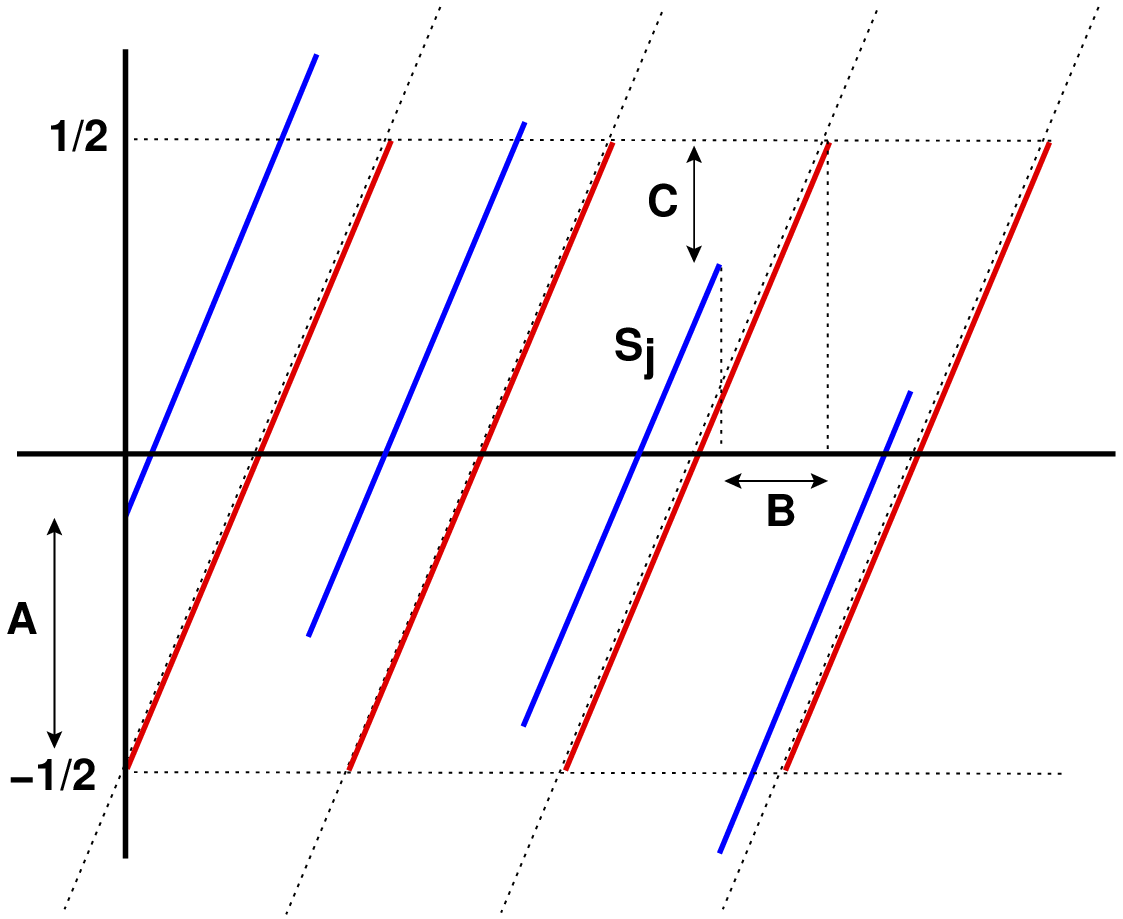}
\caption{\emph{Sketch of the branches of $S$ with $d=0$ in red and with $d>0$ in blue.
Here, $S_j$ is the third branch. The lengths of the segments indicated by $A$, $B$, and $C$,
are, respectively, $\frac{(q+1)d}{2}$, $\frac{i_+d}{q}$, and $\frac{(q+1-2i_+)d}{2}$.}}
\label{fig-branchesofG}
\end{figure}

\vskip -0.0in\noindent
{\bf Proof.} Fix $p$, $q$, and $\rho$. By definition of $\nu$, its value at any point $z$ equals
the number of solutions $x$ of $S(x)=z$ divided by $q$. Since $S$ has $q$ branches, then
$\nu(\rho,q,z)\neq 1$ implies that $\nu(\rho,q,z)\in [0,\frac{q-1}{q}]$. So all we just need to
prove the first statement of the proposition.

Our strategy is to derive explicit equations for each of the branches of $S(\rho,q,x)$ (see Figure
\ref{fig-branchesofG}). Let $\rho= \frac{p+d}{q}$. For ease of
expressing the ideas, we think of each of the branches as including \emph{both} endpoints.
\bse
S(x)= qx+\dfrac{(q+1)(p+d)-q}{2} - \sum_{i=1}^{q}\,\left\lfloor x+ i\left(\frac{p+d}{q}\right)\right\rfloor \,.
\label{eq:Gx-with-d}
\ese
Now set $x=0$ and use Proposition \ref{prop:funnysum} (i) to evaluate $S(0)$  for $d\in
\left(\frac{-1}{q-1},\frac{1}{q-1}\right)$.
\bse
S(0)=\frac{(q+1)d-1}{2}-\lfloor d \rfloor = \frac{(q+1)d-\textrm{sgn}(d)}{2} \,,
\label{eq:G0-with-d}
\ese
with the proviso that in this case we understand $\textrm{sgn}(0)=+1$. For ease of argument we
restrict to $d\geq 0$. The case $d<0$ is almost identical.

Clearly, the graph of each of the $q$ branches of $S$ is a affine segment with slope $q$. We label them from 
left to right by $S_{0}$, $S_{1}$ ,..., $S_{q-1}$. So the first $S_{0}$ is a segment of the line
$qx+ \frac{(q+1)d-1}{2}$. The discontinuity causes a reduction by one, and so $S_{1}$ lies on
the line $qx+ \frac{(q+1)d-1}{2}-1$, and so on. Thus
\bsenn
S_{j}(x)= qx+\frac{(q+1)d-1}{2}-j \,.
\esenn
The locations of the discontinuities are given by $\{-i\rho\}_{i=1}^q=\left\{-i(p+d)/q\right\}_{i=1}^q$.
Thus as $d$ increases from 0, the $i$th discontinuity moves to the \emph{left} by $\frac{id}{q}$, where
$i\in\{1,\cdots,q\}$. The left- and right discontinuities of $S_{j}$ undergo different shifts.
When $d=0$, the domain of the $j$th branch $S_j$ is $[\frac{j}{q},\frac{j+1}{q}]$. Thus the values
of $i$ are $i_-$ at the left discontinuity of $S_j$ and $i_+$ at the right one, where
\bsenn
i_+ \frac pq \mod 1= \frac jq \quad \logand \quad i_+ \frac pq \mod 1 =\frac{j+1}{q} \,.
\esenn
Therefore the domain of the $j$th branch with $d>0$ is $\left[\frac{j}{q}-\frac{i_-d}{q},
\frac{j+1}{q}-\frac{i_+d}{q}\right]$. Substituting these values into the formula for $S_j$,
\bse
S_{j}(\textrm{left endpoint}) = \frac{(q+1-2i_-)d-1}{2} \quad \logand \quad
S_{j}(\textrm{right endpoint}) = \frac{(q+1-2i_+)d+1}{2}  \,.
\label{eq:jth-branch}
\ese
Thus for $i_-$ and $i_+$ in $\{1,\cdots, q\}$, the condition on $d$ guarantees that the value at
the left endpoint of $S_{j}$ is less than 0, while the value at the right endpoint is greater than 0.
More precisely, for $d>0$, the supremum of the images of these left endpoints equals
$\frac{(q-1)d-1}{2}$ while the infimum of the images of the right endpoints equals $\frac{(1-q)d+1}{2}$,
which is what we wanted to prove.  \QED

\vskip -0.2in\noindent
\begin{prop} Let $q\geq 2$ be a positive integer, $|d|<1/(q-1)$, and $\rho, \rho'$ satisfy
\[ \rho = \frac{p+d}{q}, \qquad \rho'=\frac{p'+d}{q}\] where both $p$ and $p'$ are relatively prime to $q$. Then for all $z \in \R$, $\nu(\rho,q,z)=\nu(\rho',q,z)$.
\label{prop:reducedresidue}
\end{prop}

\vskip -0.0in\noindent
{\bf Proof.}
In light of Definition \eqref{def:density} we need only verify that the number of inverse images $x$ of
$S(\rho,q,x)=z$ equals the number of inverse images of $S(\rho',q,x)=z$. But in the proof of
Proposition \ref{prop:nu=1}, the collection of values of the Birkhoff sum at the endpoints of the
branches in \eqref{eq:jth-branch} does not depend on $p$ (as long as it is coprime to $q$). The
values at the left endpoints of the branches are permuted for $S(\rho,q,x)$ and $S(\rho',q,x)$ and
the same holds for the values at the right endpoints (though permuted differently).
This proves the proposition.
\QED

\vskip -0.2in\noindent
\begin{defn} We say that a function $v$ is a trapezoid of step $n$, if there is a sequence of positive numbers $0<t_1<t_2<\dots<t_n$, such that (see Figure \ref{fig:construction-bhoffmeas}):\\[-0.8cm]
\begin{enumerate} \itemsep -0.07cm
\item[a.] The image of $v$ is $\{t_i, 1\le i\le n\}$, and $v^{-1}{\{t_n\}}$ is a single interval (open or half open),
\item[b.] For $1\le i\le n-1$, the set $v^{-1}(t_i)$ is a union of two intervals and these $2(n-1)$ intervals are of equal length.
\item[c.] Left endpoints of $\inf v^{-1}(t_i)$ form an increasing sequence:
$\inf v^{-1}(t_1)< \inf v^{-1}(t_2)<\dots <\inf v^{-1}(t_n)$. The corresponding right endpoints
$\sup v^{-1}(t_i)$ form a decreasing sequence.
\end{enumerate}
\label{def:trapezoid}
\end{defn}

\begin{figure}[pbth]
\center
\includegraphics[width=3.5in]{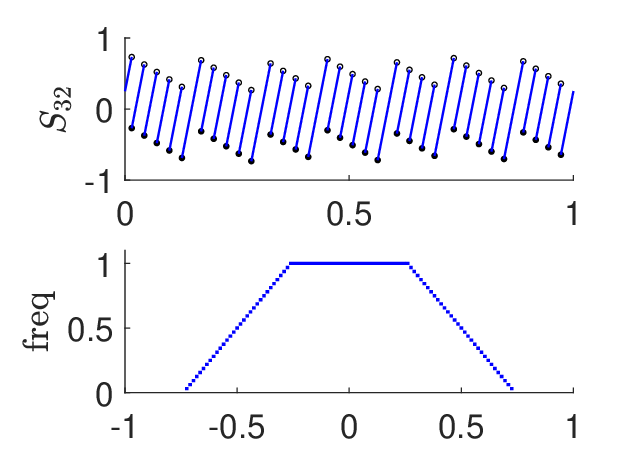}
\includegraphics[width=3.5in]{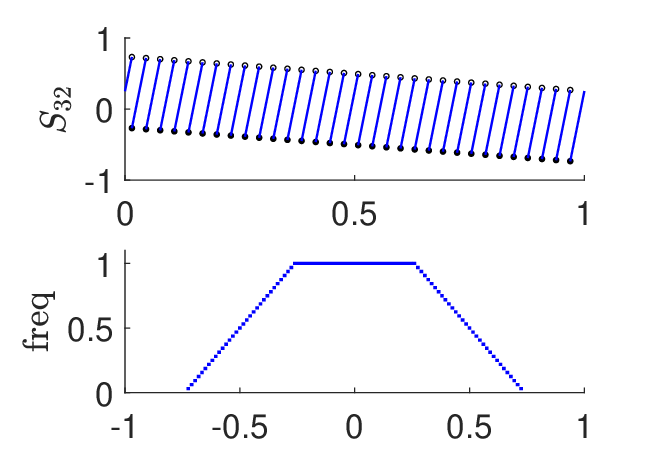}
\caption{\emph{On the left we construct $\nu(e-2,32,z)$. Knowing that $23/32$ is an odd
approximant of $e-2$, we then construct $\nu(e-2-22/32, 32,z)$ on the right. }}
\label{fig:reducedresidue}
\end{figure}

\vskip -0.0in\noindent
\begin{theo} If $p_n/q_n$ is a continued fraction convergent of $\rho$, then the graph of 
$z\rightarrow \nu(\rho,q_n,z)$ is an isosceles trapezoid of step $q_n$.
\label{thm:trapezoid}
\end{theo}

\vskip -0.0in\noindent
{\bf Proof.} First let $p/q:=p_{2n+1}/q_{2n+1}$ be an \emph{odd} approximant for $\rho$, so that
$d<0$, and set $\rho':=\frac{1+d}{q}$. By Proposition \ref{prop:reducedresidue},
$\nu(\rho',q,z)=\nu(\rho,q,z)$. We now have $q\rho'-1=d<0$.
Therefore the continued fraction expansion of $\rho'$ starts with $q$: $[q,\cdots]$.
So if we use the conventions of Definition \ref{def:contdfr} and \ref{def:dn} for $\rho'$, we now have
\bsenn
p_0/q_0=0/1 \quad \logand \quad p_1/q_1=1/q  \quad \logand \quad d_0=\rho' \quad \logand \quad d_1=q\rho'-1 \,.
\esenn
Using the special case of the \emph{three gaps theorem} \cite{book}, Proposition B.3, we see that the
spacings between successive points $\{i\rho\}_{i=1}^{q_1}$ come in two sizes, namely one ``large" gap of size
$d_0 + |d_1|$ and $q-1$ gaps of size $|d_1|$, where the large gap (on the circle) is comprised of the
union of $[0,|d_1|)$ and $[1-d_0,1)$. This implies that the upper tips of the branches of $S(\rho',q,x)$
are equidistant points lying on a straight line $\ell$, while the lower limits of the branches
lie on a translate of $\ell$ by one. Since $\nu$ is the pushforward of the Lebesgue measure by $S$,
this establishes that $\nu$ is an isosceles trapezoid when $p/q$ is an odd approximant. The transition
from $p$ to $p'=1$ is illustrated in Figure \ref{fig:reducedresidue}.

If $p/q:=p_{2n}/q_{2n}$ is an \emph{even} approximant for $\rho$, by Proposition \ref{prop:symmetry}, we
may equivalently consider $\nu(1-\rho,q,z)$. It is now easy to show that $q$ is the denominator
of an odd approximant. \QED

\begin{figure}[pbth]
\center
\includegraphics[width=4.2in]{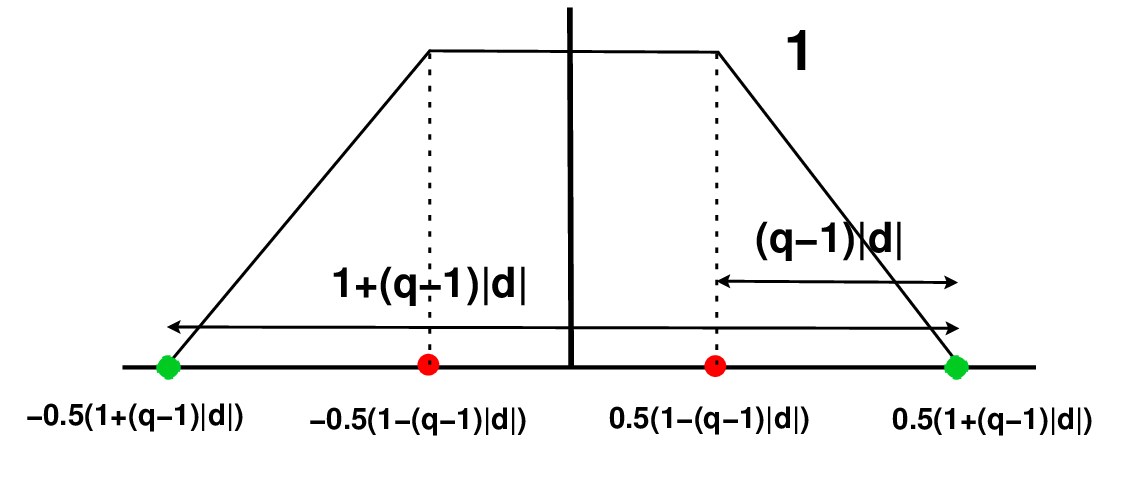}
\caption{\emph{Illustration of the shape of $z\rightarrow \nu(\rho,q_n,z)$ (Corollary \ref{cory:trapezoid}). }}
\label{fig:trapezoid}
\end{figure}

\vskip -0.0in\noindent
\begin{cory}[The Trapezoid Theorem] If $p_n/q_n$ is a continued fraction convergent of $\rho$,
then $z\rightarrow \nu(\rho,q_n,z)$ is an isosceles trapezoid of step $q_n$ (a trapezoid that is symmetric under
$z\rightarrow -z$) which equals 1 on an interval of length $1-(q_n-1)|d|$ and with support of
length $1+(q_n-1)|d|$ (see Figure\ref{fig:trapezoid}).
\label{cory:trapezoid}
\end{cory}

\vskip -0.0in\noindent
{\bf Proof.} The first statement follows from Theorem \ref{thm:trapezoid}, the second from Proposition
\ref{prop:nu=1}, and the third from the fact that the integral of $\nu(\rho,q_n,z)$ equals one. \QED

\vskip -0.1in\noindent
{\bf Remark.} Together with Corollary \ref{cor:Cn-is-length-of-support of nu}, this last result implies
that $q_nD_{q_n}(\overline{x}_\rho)$ equals $1+(q_n-1)|d_n|$ where $q_n$ and $d_n$ are as in Definition
\ref{def:contdfr} and \ref{def:dn} and $\overline{x}_\rho=\{i\rho\}_{i=1}^\infty$.

\begin{centering}
\section{Appendix A: The Relation Between Discrepancy and Birkhoff Sums}
\label{chap:Clump+Bhoff}
\end{centering}
\setcounter{figure}{0} \setcounter{equation}{0}

\vskip 0.0in\noindent
In this section, we state a result from \cite{Ramshaw} and give a new, much shorter proof. We
will repeatedly use the following two observations. The maximum and the minimum of the range of
$S(\rho,n,x)$ occur at the discontinuities located at $\{-i\rho\}$. Furthermore,
\bse
\forall\;x\not\in \Z:\;\;\{-x\}=1-\{x\} \quad \logand\quad  \forall\;x\in \Z:\;\;\{-x\}=\{x\}=0 \,.
\label{eq:frac-neg-x}
\ese

\vskip -0.0in\noindent
\begin{lem} Let $\rho$ be irrational and $0<\ell, k\leq n$. Then
\bsenn
\sum_{i=1}^n\left(\{(i-k)\rho\}-\frac 12\right)+\sum_{i=1}^n\left(\{(i-\ell)\rho\}-\frac 12\right)=-1
\quad \Longleftrightarrow \quad k+\ell=n+1 \,.
\esenn
\label{lem:addingsequences}
\end{lem}

\vskip -0.4in\noindent
\begin{figure}[pbth]
\center
\includegraphics[width=4.5in]{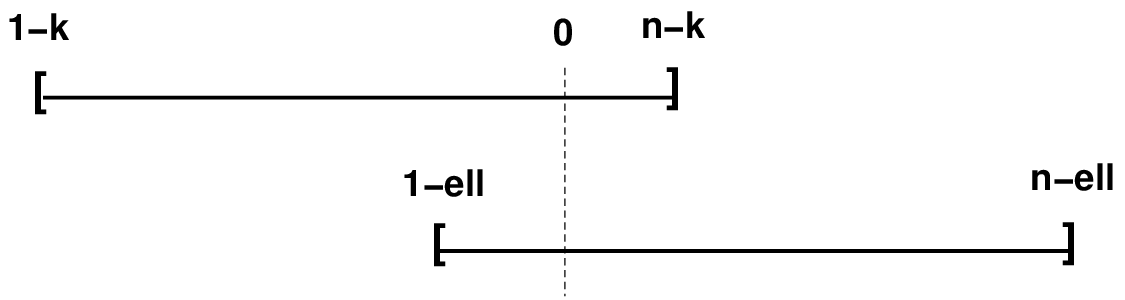}
\caption{\emph{Adding the two sequences of Lemma \ref{lem:addingsequences}. In order for all
the irrational terms to cancel, the picture has to be symmetric around zero.}}
\label{fig:addingsequences}
\end{figure}

\vskip -0.0in\noindent
{\bf Proof.} We may assume without loss of generality that $\ell\leq k$. Relabel the indices in the sums as follows (see Figure \ref{fig:addingsequences}).
\bsenn
\sum_{i=1}^n\left(\{(i-k)\rho\}-\frac 12\right)+\sum_{i=1}^n\left(\{(i-\ell)\rho\}-\frac 12\right) =
\sum_{j=1-k}^{n-k}\left(\{j\rho\}-\frac 12\right)+\sum_{j=1-\ell}^{n-\ell}\left(\{j\rho\}-\frac 12\right) \,.
\esenn
From \eqref{eq:frac-neg-x}, we see that the two terms with $j=0$ together give a contribution -1.
Since $\rho$ is
irrational, all terms with $j\neq 0$ satisfy $\{-j\rho\}-\frac 12=\frac 12-\{j\rho\}$. We use
irrationality of $\rho$ again by asserting that the $j\neq 0$ contributions cancel if and only if
all those terms occur in pairs $j$ and $-j$. From Figure \ref{fig:addingsequences}, we see that
this happens if and only if $-(1-k)=n-\ell$ and $1-\ell=n-k$. In turn, these are equivalent to
$k+\ell=n+1$.
\QED

\vskip -0.2in\noindent
\begin{prop} The minimum of $S(x) $ is achieved at $x=\{-k\rho\}$ for some $k\in\{1,\cdots,n\}$,
while its maximum is achieved at $x=\{-\ell\rho\}$ with $k+\ell=n+1$.
\label{prop:minandmax of S}
\end{prop}

\vskip 0.0in\noindent
{\bf Proof.} Clearly, the extrema of $S$ are located at the discontinuities. Therefore the global minimum
of $S$ is $\sum_{i=1}^n\left(\{(i-k)\rho\}-\frac 12\right)$ for some $k$ while the global maximum equals
a local minimum plus one, and so must be equal to
$\sum_{i=1}^n\left(\{(i-\ell)\rho\}-\frac 12\right)+1$ for some $\ell$. Since the range of $S$ is
symmetric, we have
\bsenn
\sum_{i=1}^n\left(\{(i-k)\rho\}-\frac 12\right)=-\left(\sum_{i=1}^n\left(\{(i-\ell)\rho\}-\frac 12\right)+1\right)\,.
\esenn
By Lemma \ref{lem:addingsequences}, $\ell$ must be equal to $n+1-k$. \QED

\vskip -0.2in\noindent
\begin{theo}[\cite{Ramshaw}] For $\overline{x}=\{i\rho\}_{i=1}^\infty$
\bsenn
nD_n(\overline{x})=1+2\,\max_{0<m < n}\left(S(\rho,m,0)-S(\rho,n-1-m,0)\right) \,.
\esenn
\label{thm:subtractingsequences}
\end{theo}

\vskip -0.2in\noindent
{\bf Proof.} In view of Theorem \ref{thm:clumps=rangeS} and Proposition \ref{prop:minandmax of S},
we need to compute
\benn
nD_n(\overline{x}) &=& 1+ \max_{k+\ell=n+1}\sum_{i=1}^n
\left(\{(i-k)\rho\}-\frac 12\right)-\sum_{i=1}^n\left(\{(i-\ell)\rho\}-\frac 12\right) \\
&=& 1+\max_{k+\ell=n+1}\sum_{j=1-k}^{n-k}\left(\{j\rho\}-
\frac 12\right)-\sum_{j=1-\ell}^{n-\ell}\left(\{j\rho\}-\frac 12\right) \,,
\eenn
after relabeling as in Lemma \ref{lem:addingsequences}. We now need to \emph{subtract} the two sequences
that we added in Lemma \ref{lem:addingsequences}. This time, the middle part
$j\in \{1-\ell,\cdots, n-k\}$ completely cancels and we end up with
\bsenn
\cdots = 1+ \max_{k+\ell=n+1}\sum_{j=1-k}^{-\ell}
\left(\{j\rho\}-\frac 12\right)-\sum_{j=n-k+1}^{n-\ell}\left(\{j\rho\}-\frac 12\right)
=1+ 2\max_{k+\ell=n+1}[S(\rho,n-\ell,0)-S(\rho,n-k+1,0)] \,.
\esenn
Substituting $m:=k-1$ gives the final result. \QED

\vskip 0.0in\noindent
\begin{centering}
\section{Appendix B: The Structure of Birkhoff Sums}
\label{chap:Bhoffrecursion}
\end{centering}
\setcounter{figure}{0} \setcounter{equation}{0}

\vskip -0.0in\noindent
The aim of this section is to state a result from \cite{DT} with a modified and shorter proof. We also
summarize a remarkable result from \cite{Ramshaw} that can be obtained from it. From here on out, we will
fix an irrational rotation number $\rho$ and the initial condition $x=0$. To simplify notation, we
write $S(i)$ for $S(\rho,i,0)$ in the remainder.

\vskip -0.0in\noindent
\begin{lem} For irrational $\rho$:
\bsenn
S(q_n)=\;\frac 12 ((q_n+1)d_n+(-1)^{n+1}) \,.
\esenn
\label{lem:contdfr2}
\end{lem}

\vskip -0.2in\noindent
{\bf Proof.} The equality follows from Proposition \ref{prop:funnysum} (ii) and the definitions of $q_n$ and $d_n$.
\QED

\vskip -0.1in\noindent
The following proposition is the main vehicle to extract the recursive structure of the sequence $S(\rho, i, 0)$.

\vskip -0.0in\noindent
\begin{prop} For irrational $\rho$:
\bsenn
\forall\,\;0\leq k<q_{n+1},\;\;\;\; S(q_n+k)= S(q_{n})+ S(k)+ kd_{n} \,.
\esenn
\label{prop:recursionS}
\end{prop}

\vskip -0.2in\noindent
{\bf Proof.} First place $\{i\rho\}$ for $i$ from 1 to $q_n$ on the unit interval to compute the
Birkhoff sum $S(\rho, q_n, 0)$.
Let $k < q_{n+1}$, and consider the set of points given by $\{i\rho\}_{i=q_n+1}^{q_n+k}$.
This is identical to the set of points $\{i\rho\}_{i=1}^k$ translated by $d_n=q_n\rho-p_n$
\emph{if and only if} each of the latter points has the property that $[i\rho,i\rho+d_n]$ (or
$[i\rho+d_n,i\rho]$ if $d_n$ is negative) does \emph{not} contain an integer. In turn, the closest
return theorem for continued fractions \cite{book}{ Chapter 6} guarantees that this is the case
if and only if $k < q_{n+1}$.
\QED

\noindent
Let $\rho=[a_1,a_2,\cdots]$ be irrational with continued fraction denominators $\{q_k\}$.
Let $L$ be an arbitrary positive integer. Its \emph{Ostrowski expansion} expresses $L$ as a finite sum:
\bsenn
L=\sum_{i=0}^nb_{i}q_i \,.
\esenn
with each `digit' $b_i$ a non-negative integer, subject to these rules: \\
(i) $0\leq b_0\leq a_1-1$ and for $i\ge 1$, $0\leq b_i\leq a_{i+1}$, and \\
(ii) if $b_{i}=a_{i+1}$, then $b_{i-1}=0$.\\
It is well known that this expansion is unique.
Given a sequence of digits $\{b_i\}_{i=0,\cdots,n}$, denote for $k\leq n$:
\bsenn
L_k:=\sum_{i=0}^kb_{i}q_i <q_{k+1} \,.
\esenn

\vskip -0.2in\noindent
\begin{prop} Let $\rho$ be irrational with continued fraction denominators $\{q_k\}$. Given the
Ostrowski representation of a positive integer as $L_n=\sum_{i=0}^nb_iq_i$, then
\benn
S(L_n) &=& \sum_{k=0}^n\;b_k\left[S(q_k)+\frac 12 (b_k-1)q_kd_k+L_{k-1}d_k\right] \\
 &=& \sum_{k=0}^n\;\frac{b_k}{2}\left[\big(b_kq_k+2L_{k-1}+1\big)d_k+(-1)^{k+1}\right]
\eenn
\label{prop:repeatedrecursionS}
\end{prop}

\vskip -0.1in\noindent
{\bf Proof.} Apply Proposition \ref{prop:recursionS} with $i=(b_k-1)q_k+L_{k-1}$:
\bsenn
S(L_k)=S(q_k+(b_k-1)q_k+L_{k-1})=S(q_k)+\left((b_k-1)q_k+L_{k-1}\right)d_k+S((b_k-1)q_k+L_{k-1}) \,.
\esenn
Repeated application from $S((b_k-1)q_k+L_{k-1})$, $S((b_k-2)q_k+L_{k-1})$, to
$S(q_k+L_{k-1})$, gives
\benn
S(L_k)-S(L_{k-1}) &=& b_k\left[S(q_k)+\left(\frac 12 (b_k-1)q_k+L_{k-1}\right)d_k\right] \\
 &=& \frac{b_k}{2}\left[\big(b_kq_k+2L_{k-1}+1\big)d_k+(-1)^{k+1}\right]
\eenn
The last equality follows by substituting Lemma \ref{lem:contdfr2}.  Summing over $k$
from 0 to $n$ gives the final result. \QED

\vskip -0.1in\noindent
{\bf Remark.} Note that this recursion reveals a possible mechanism to select irrationals and
sequences $L_n$ for which $S(L_n)$ can grow with almost linear growth in $n$ by choosing $\limsup
\frac{b_{n+1}}{L_n}$ to be
positive.

\vskip 0.1in
It is important to determine the influence of one coefficient $b_m$ on the outcomes of
all $S(L_n)$ for all $n\geq m$. For $m=n$, $b_m$ appears explicitly in the last
term of the sum in Theorem \ref{prop:repeatedrecursionS}. However, when $m<n$, $b_m$ also
appears `hidden' in the terms $L_{k-1}$ with $m\leq k-1< n$.
The following Theorem, equivalent to Proposition 1.71 of \cite{DT}, makes this dependence explicit
and is the most important tool in the computation of Birkhoff sums.

\begin{figure}[!ht]
\centering
\includegraphics[height=2.5in,width=3.4in]{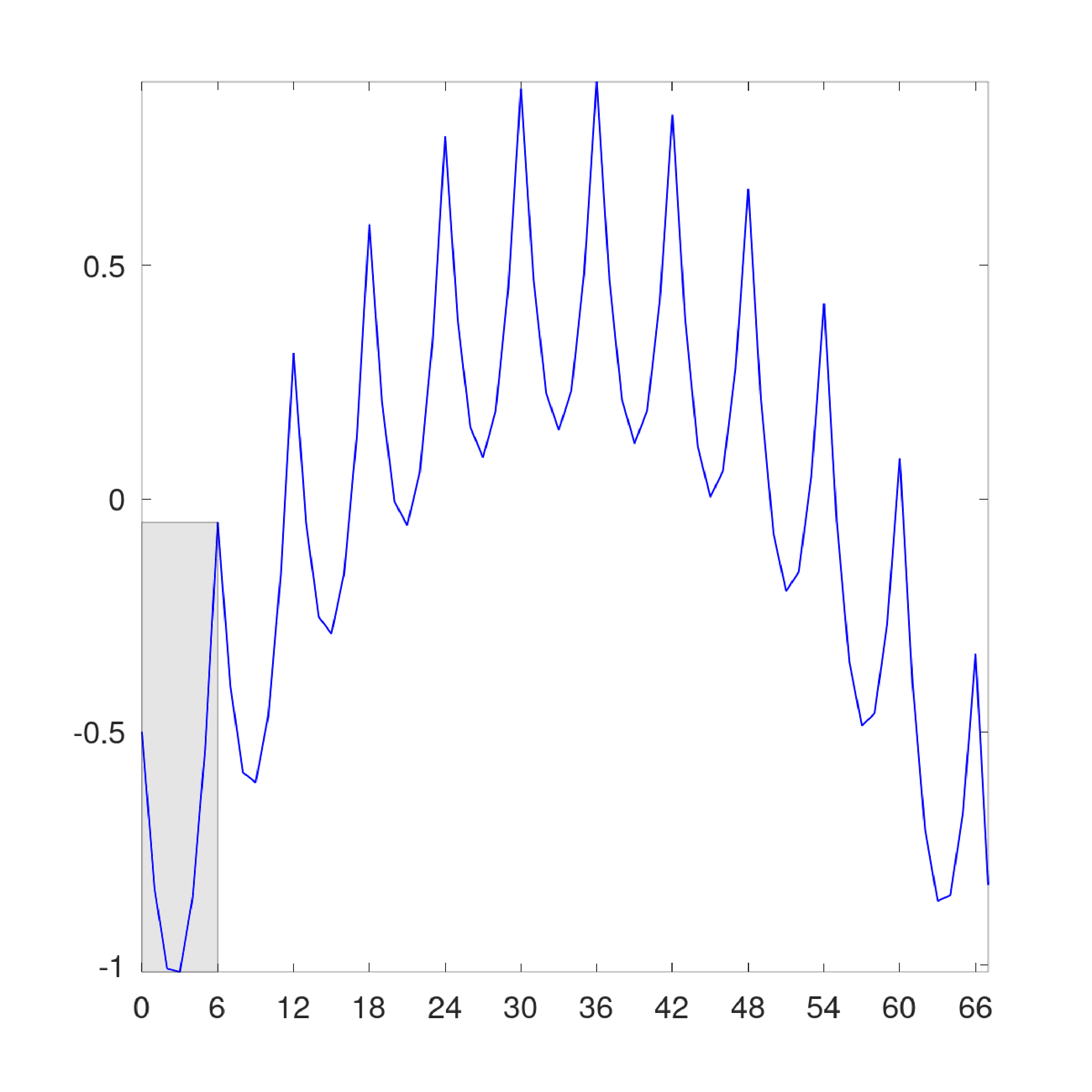}
\includegraphics[height=2.5in,width=3.4in]{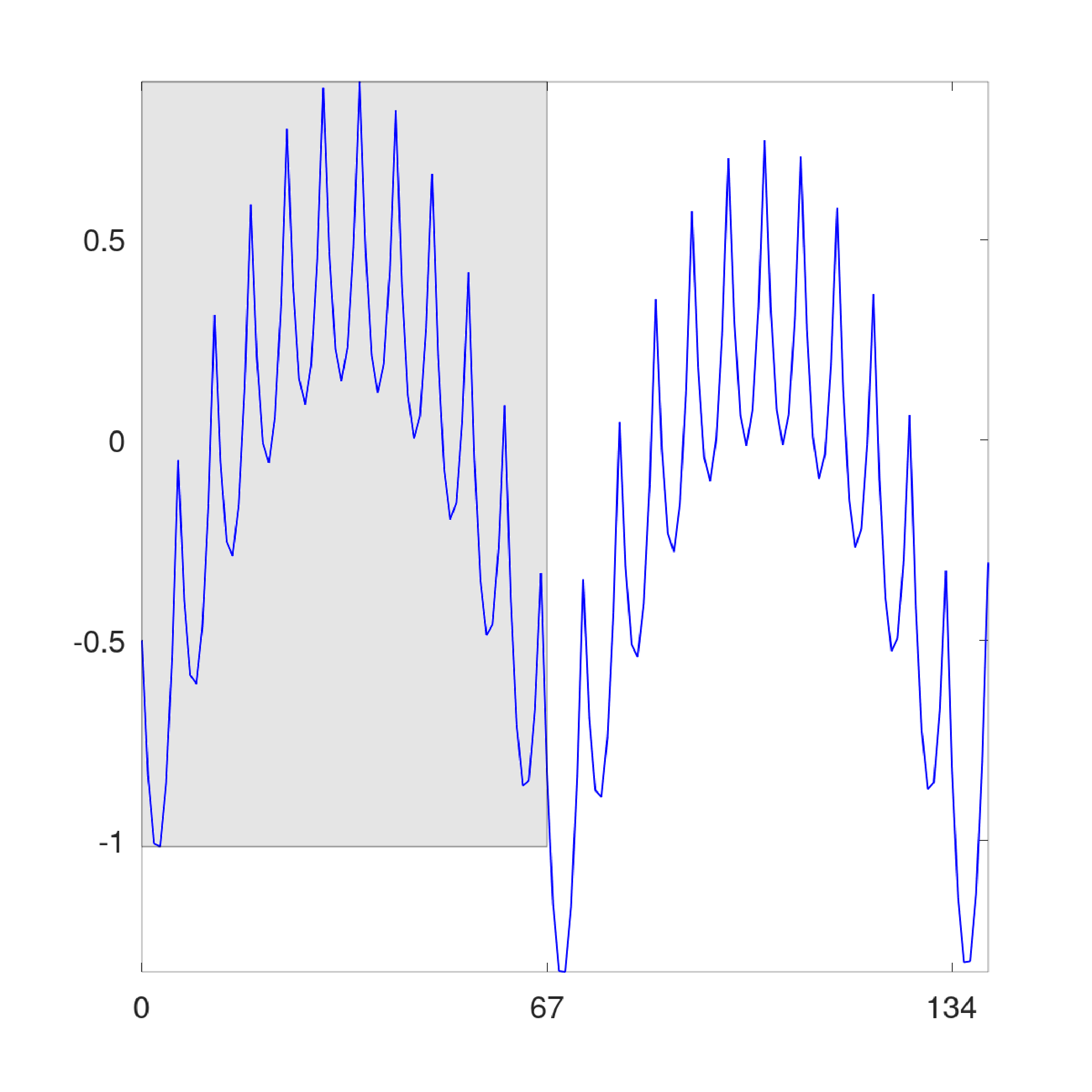}
\includegraphics[height=2.5in,width=3.4in]{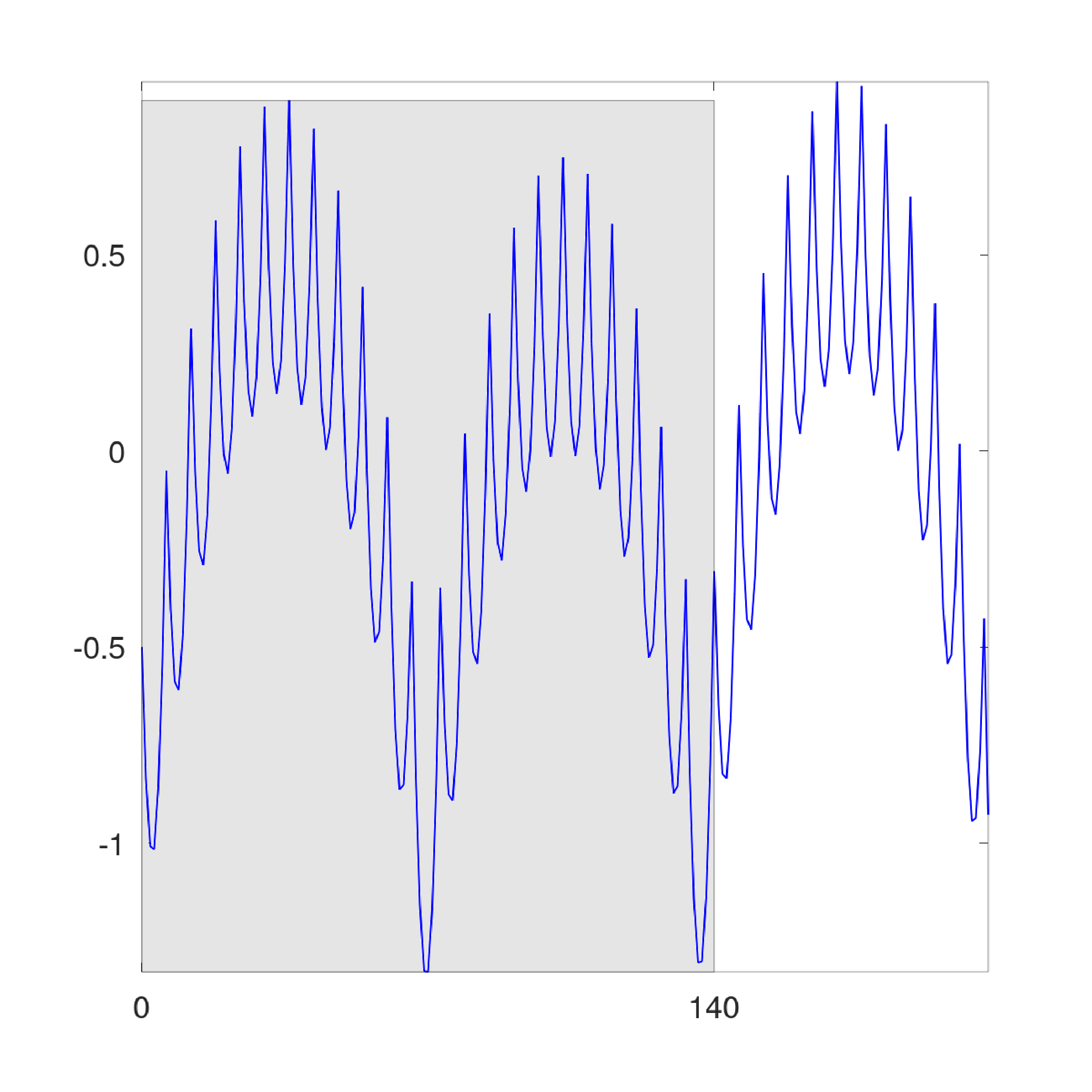}
\includegraphics[height=2.5in,width=3.4in]{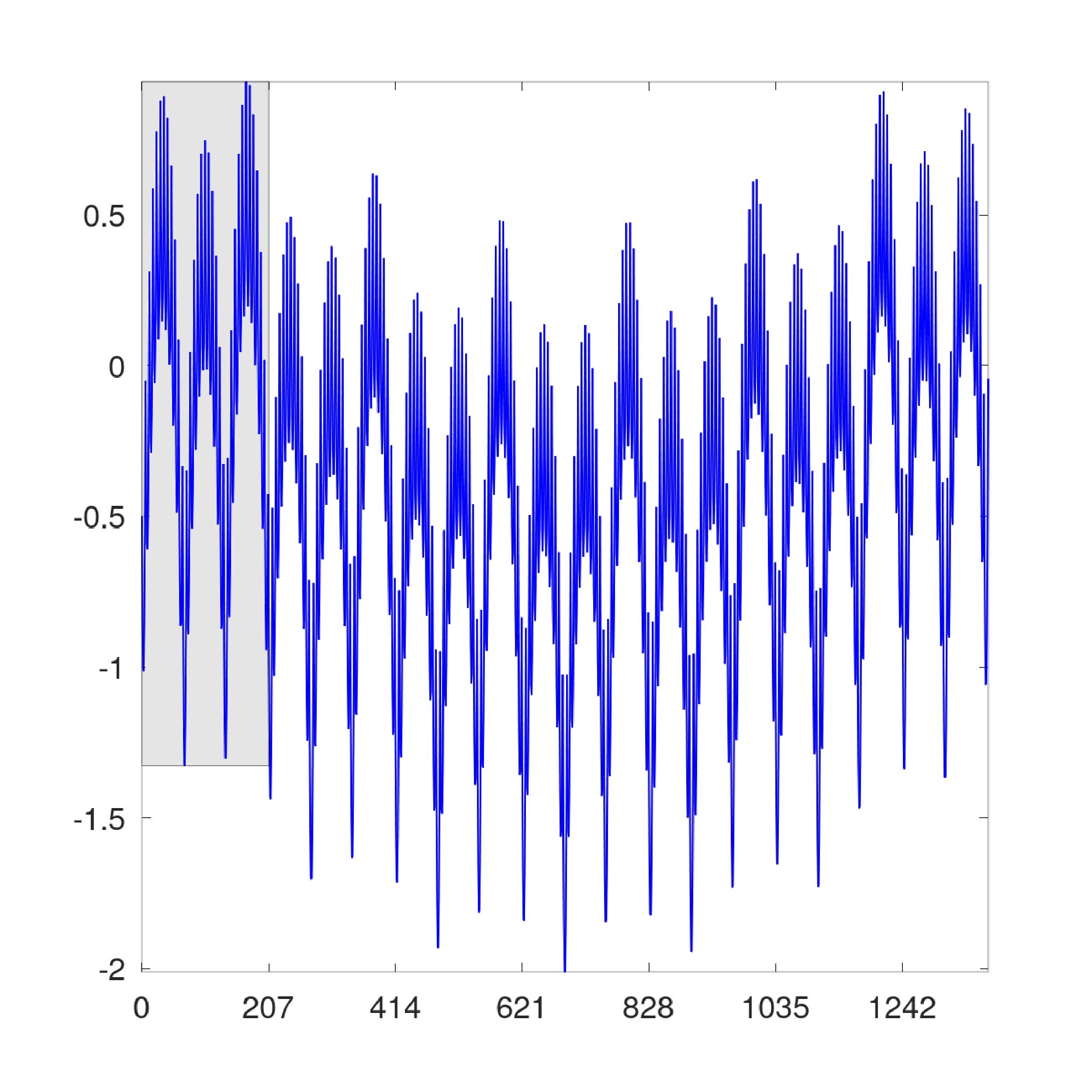}
\caption{\emph{Birkhoff sums for the repeated continued fraction $\rho=[\overline{6,11,2,1}]$.
The parabola-like structure formed by the $a_n$ affine translations is apparent when $a_n$ is large
(e.g. for $q_2$ and $q_5$) but less clear when $a_n$ is small (e.g $q_3$ and $q_4$).}}
\label{fig:fractabolae}
\end{figure}

\vskip -0.1in\noindent
\begin{theo} Let $\rho$ be irrational with continued fraction denominators $\{q_k\}$. Given the
Ostrowski expansion $L_n=\sum_{i=0}^nb_iq_i$. Suppose for some $m\leq n$, we vary \underbar{only}
$b_m$ and keep $b_j$ for $j\neq m$ constant. Then
\bsenn
S(L_n)=c+ B(m,n,b_m) \,,
\esenn
where $c$ is a function of the $b_j$, $j\neq m$ and $B(m,n,b_m)$ satisfies
\begin{equation}
B(m,n,b_m)= \frac{b_m^2}{2}\;q_md_m+b_m\;\left(\frac{(-1)^{m+1}}{2}+L_{m-1}d_m+
\frac{d_m}{2}+q_{m}\sum_{j=m+1}^nb_{j}d_j\right) \,.
\label{eq:B}
\end{equation}
\label{thm:dependence-bk}
\end{theo}

\vskip -0.2in\noindent
{\bf Proof.} The main change in the second equality of Proposition \ref{prop:repeatedrecursionS}
is the computation of the coefficient of the variable $b_m$ in the quadratic form $\sum_{k=0}^n b_k L_{k-1} d_k$ which equals $$b_m L_{m-1} d_m+b_m q_m \sum_{j=m+1}^n b_j d_j$$.
The rest is a minor rewrite to bring out the quadratic in  terms of the variable $b_m$. \QED

\vskip -0.1in\noindent
The expression $B(m,n,b_m)$ is quadratic in $b_m$ and is used to represent the influence
of the $m$th digit $b_m$ on the value of $S(L_m)$. This quadratic form causes the Birkhoff sums to have a
combined fractal and parabolic character, resulting in shapes we call `fractabolae', as shown in Figure
\ref{fig:fractabolae}. In the figure, $\rho$ has repeated continued fraction $[\overline{6,11,2,1}]$ and
the $q_i$ are (starting with $q_0$): 1, 6, 67, 140, 207, 1382, ... . The sums through the previous
$q_{n-1}$ are in a shaded rectangle. For odd values of $m$ the
coefficient $q_md_m$ of the quadratic term is negative and the quadratic is concave and has a maximum
value. For even values of $m$ the quadratic is convex with a minimum value.

\vskip -0.0in
To illustrate what can be computed with this classical result, we finish this section by giving an
example. Let $\rho_a=[a,a,a,\cdots]$ with $a\in \N$. Examples are the golden mean
$\rho_1=\frac{\sqrt{5}-1}{2}=[1,1,....]$ and the silver mean
$\rho_2=\sqrt{2}-1=[2,2,....]$. For these numbers, one can use the machinery in this section
to prove Ramshaw's result \cite{Ramshaw}, namely that for these `metallic' means,
$\limsup_{n\rightarrow \infty}\dfrac{S(\rho_a,n,0)}{\ln n}$ exists and is equal to
\bsenn
c(a)=
\begin{cases}
 \dfrac{a}{16}\dfrac{1}{\ln \rho_a^{-1}} & \logif \;a \;\; {\rm is \; even}\\[-0.3cm]\\
 \dfrac{a(a^2+3)}{16(a^2+4)}\dfrac{1}{\ln \rho_a^{-1}} & \logif \;a \;\; \rm{is\;odd} \,.
\end{cases}
\esenn
As a remarkable consequence, see Corollary \ref{cor:Cn-is-length-of-support of nu},
$\limsup_{n\rightarrow \infty}\dfrac{nD_n(\{i\rho_a\}_{i=1}^n)}{\ln n}$  exists and is equal to $4c(a)$
\cite{Ramshaw}. Nonetheless, it follows from the ergodicity of the Gauss map that for Lebesgue
almost all $\rho$,
\bsenn
\limsup_{n\rightarrow \infty}\dfrac{S(\rho,n,0)}{\ln n}=\infty \,.
\esenn
Thus, for these $\rho$, the $\limsup$ of $nD_n(\{i\rho\}_{i=1}^n)/\ln n$ also tends to infinity.

%\vskip 0.4in\noindent
%\newpage
\begin{centering}
\section{Appendix C: A Sampling of Birkhoff Measures}
\label{chap:Cbestiary}
\end{centering}
\setcounter{figure}{0} \setcounter{equation}{0}

\begin{figure}[!ht]
\centering
\includegraphics[width=2.1in]{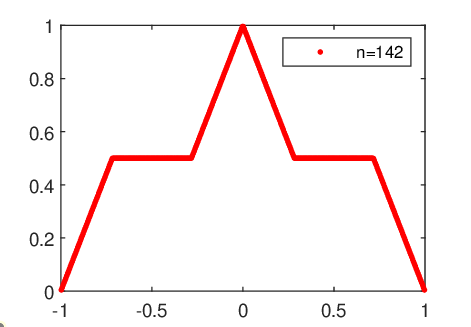}
\includegraphics[width=2.1in]{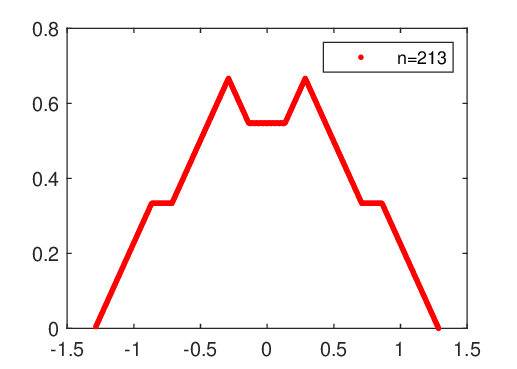}
\includegraphics[width=2.1in]{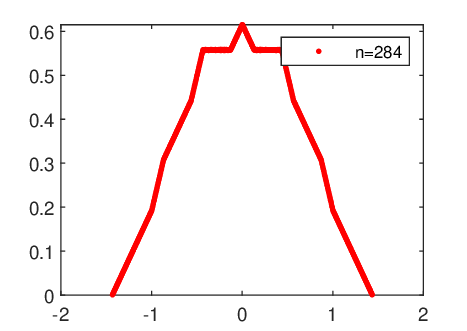}
\includegraphics[width=2.1in]{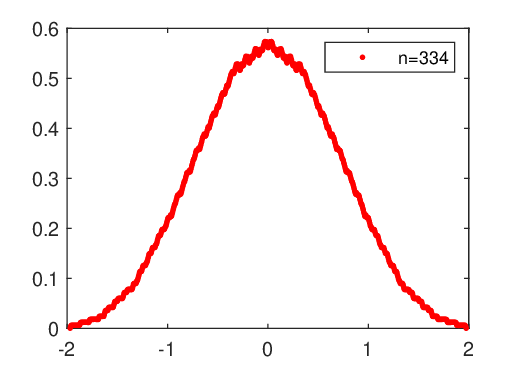}
\includegraphics[width=2.1in]{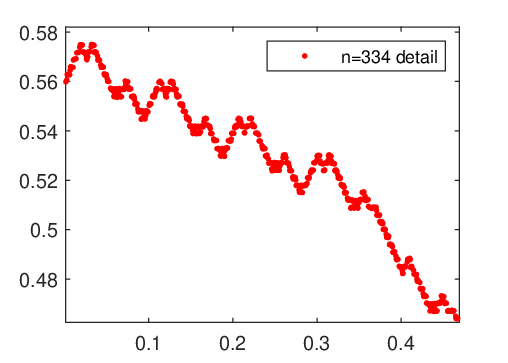}
\includegraphics[width=2.1in]{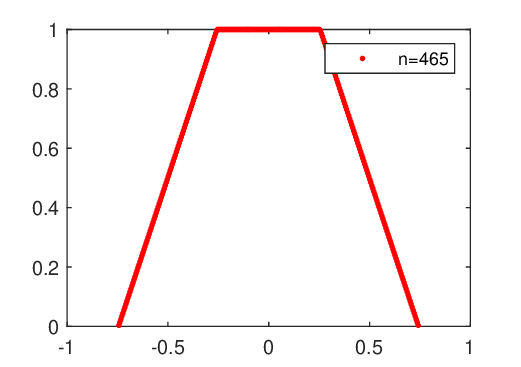}
\includegraphics[width=2.1in]{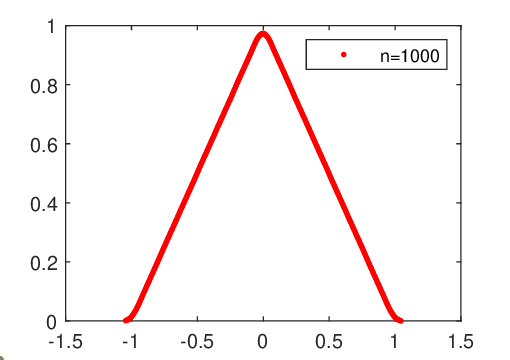}
\includegraphics[width=2.1in]{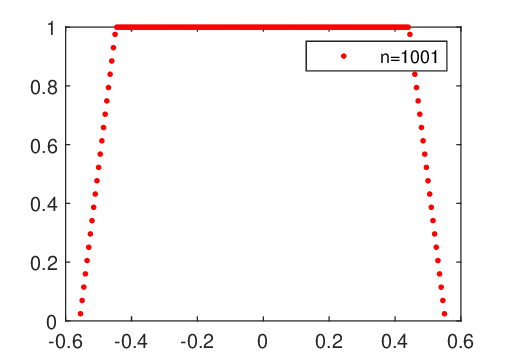}
\includegraphics[width=2.1in]{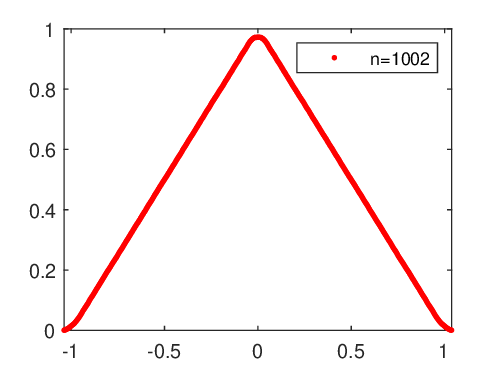}
\includegraphics[width=2.1in]{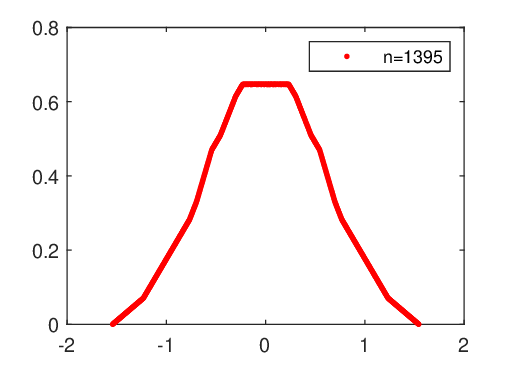}
\includegraphics[width=2.1in]{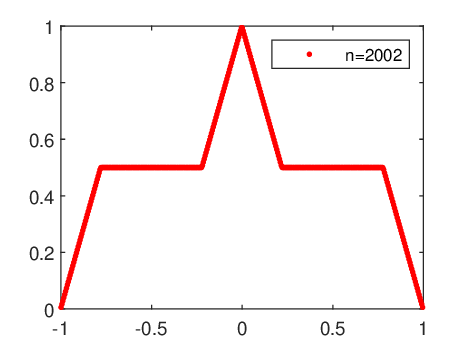}
\includegraphics[width=2.1in]{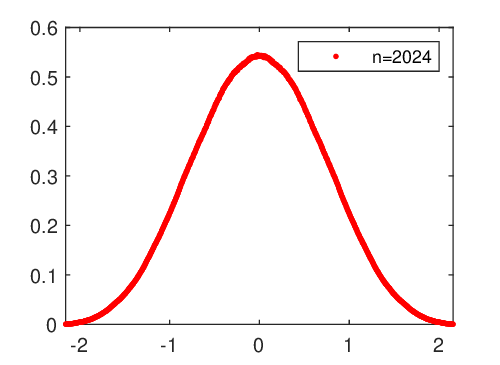}
\includegraphics[width=2.1in]{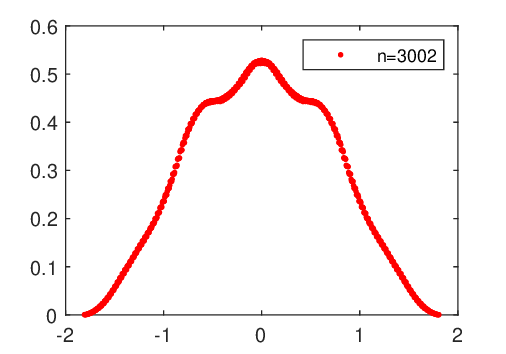}
\includegraphics[width=2.1in]{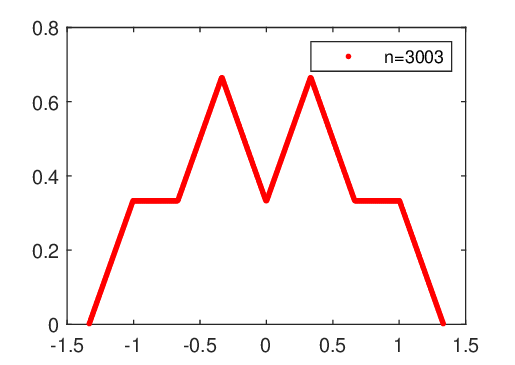}
\includegraphics[width=2.1in]{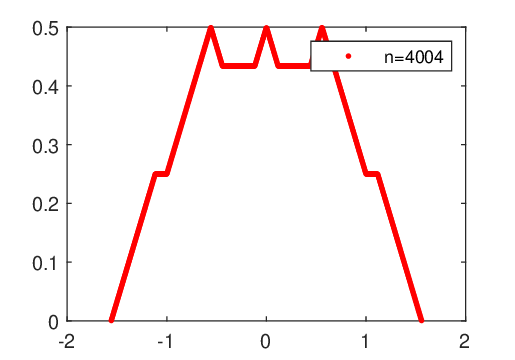}
\caption{\emph{Some Birkhoff measures for $\rho=e-2$. }}
\label{fig:measurese-2}
\end{figure}

In Figure \ref{fig:measurese-2}, we give an idea of the stunning
variety of these densities for $\rho=e-2$. Recall that $e-2=[1,2,1,1,4,1,1,6,1,1,8,\cdots]$. The
approximants are: $1, \frac{2}{3}, \frac{3}{4}, \frac{5}{7}, \frac{23}{32},\frac{28}{39}, \frac{51}{71},
\frac{334}{465}, \frac{385}{536}, \frac{719}{1001}, \frac{6137}{8544}, ...$. The middle row
shows the densities $\nu(\rho,n,z)$ when $n$ equals a continued fraction denominator (1001).
Note that even changing $n$ by $\pm 1$ can cause dramatic changes in the graph of the density.

\vspace{\fill}
\end{document}